%\special{src: 1 Y3.TEX}
\input amstex

\input epsf
\epsfverbosetrue \input amsppt.sty
\magnification 980\vsize=21 true cm \hsize=16 true cm \voffset=1.1
true cm \pageno=1 \NoRunningHeads \TagsOnRight
%\scalelinespacing{\magstep1.5}
%\special{src: 11 Y3.TEX}
\def\p{\partial}
\def\ve{\varepsilon}
\def\f{\frac}

\def\na{\nabla}
\def\la{\lambda}
\def\al{\alpha}

\def\t{\tilde}

\def\o{\omega}
\def\th{\theta}

\def\g{\gamma}
\def\G{\Gamma}
\def\si{\sigma}

\def\dl{\delta}

\def\ds{\displaystyle}

%% ------------------------------------------------------------------------

\topmatter

%% ------------------------------------------------------------------------

\title On the lifespan of and the blowup mechanism for smooth solutions to
a class of 2-D nonlinear wave equations with small initial data
\endtitle

\author DING Bingbing, Ingo WITT, and YIN Huicheng\endauthor

\address Department of Mathematics and IMS, Nanjing University,
Nanjing 210093, P.R.~China \endaddress

\email balancenjust\@yahoo.com.cn \endemail

\thanks Ding Bingbing and Yin Huicheng were supported by the NSFC
(No.~10931007, No.~11025105), by the Priority Academic Program
Development of Jiangsu Higher Education Institutions, and by the DFG
via the Sino-German project ``Analysis of PDEs and Applications.''
Some parts of this work were done when Ding Bingbing and Yin Huicheng
were visiting the Mathematical Institute of the University of
G\"{o}ttingen.\endthanks

\address Mathematical Institute, University of G\"ottingen,
Bunsenstr.~3-5, D-37073 G\"ottingen, Germany \endaddress

\thanks Ingo Witt was supported by the DFG via the
Sino-German project ``Analysis of PDEs and Applications.'' \endthanks

\email iwitt\@uni-math.gwdg.de \endemail

\address Department of Mathematics and IMS, Nanjing University,
Nanjing 210093, P.R.~China \endaddress

\email huicheng\@nju.edu.cn \endemail



%% ------------------------------------------------------------------------

\abstract
This paper is concerned with the lifespan and the blowup mechanism
for smooth solutions to the 2-D nonlinear wave equation
$\p_t^2u-\ds\sum_{i=1}^2\p_i(c_i^2(u)\p_iu)$ $=0$, where $c_i(u)\in
C^{\infty}(\Bbb R^n)$, $c_i(0)\neq 0$, and
$(c_1'(0))^2+(c_2'(0))^2\neq 0$. This equation has an interesting
physics background as it arises from the pressure-gradient model in
compressible fluid dynamics and also in nonlinear variational wave
equations. Under the initial condition $(u(0,x), \p_tu(0,x))=(\ve
u_0(x), \ve u_1(x))$ with $u_0(x), u_1(x)\in C_0^{\infty}(\Bbb R^2)$,
and $\ve>0$ is small, we will show that the classical solution
$u(t,x)$ stops to be smooth at some finite time $T_{\ve}$. Moreover,
blowup occurs due to the formation of a singularity of the first-order
derivatives $\na_{t,x}u(t,x)$, while $u(t,x)$ itself is continuous up
to the blowup time~$T_{\ve}$.
\endabstract

%% ------------------------------------------------------------------------

\endtopmatter

%% ------------------------------------------------------------------------

\document

%% ------------------------------------------------------------------------
\head  \S 1. Introduction and main result \endhead

In this paper, we are concerned with the lifespan $T_{\ve}$ of and the
blowup mechanism for classical solutions to the 2-D nonlinear
wave equation
$$
\cases &\p_t^2u-\ds\sum_{i=1}^2\p_i(c_i^2(u)\p_iu)=0,\\
&(u(0,x), \p_tu(0,x))=(\ve u_0(x), \ve u_1(x)),
\endcases\tag1.1
$$
with small initial data, where $c_i(u)\in C^{\infty}(\Bbb R^n)$,
$c_i(0)\neq0$, and $(c_1'(0))^2+(c_2'(0))^2\neq0$. In addition,
$u_0(x), u_1(x)\in C_0^{\infty}(\Bbb R^2)$ and $\ve>0$ is sufficiently
small.

In case $c_1(u)=c_2(u)=e^{u/2}$, Eq.~(1.1) has a background in physics
as it arises from the pressure-gradient model in compressible Euler
systems and is rather analogous to the 2-D variational wave equation
$\p_t^2u-\ds\sum_{i=1}^2c(u)\p_i(c(u)\p_iu)=0$ (for the physics
backgrounds of the variational wave equation and its mathematical
treatment, see [2, 10, 15, 24] and the references therein).

Here is a derivation of the pressure-gradient model with small initial
data: As pointed out in [1, 30-31], the pressure-gradient system is a
simplified version of the compressible Euler equations, which arises
from splitting the compressible Euler system (i.e., the inertia terms
$\operatorname{div}(\rho U)$, $\operatorname{div}(\rho U\otimes U)$
and the pressure $p$ are considered separately). It has the form
$$
\cases
&\p_t\rho=0,\\
&\p_t(\rho U)+\na p=0,\\
&\p_t(\rho E)+\operatorname{div}(p U)=0,
\endcases\tag1.2
$$
where $\rho$ is density, $U=(u_1, u_2)$ is velocity, $p$ is pressure,
$E=\ds\f{1}{2}\,|U|^2 +\ds\f{1}{\g -1}\f{p}{\rho}$ is energy, and $\g$
is the adiabatic exponent with $1<\g<3$.

For simplicity, as in [19-20, 25], we assume $\rho\equiv 1$ in
(1.2). In this case, (1.2) becomes
$$
\cases
&\p_t U+\na p=0,\\
&\p_t E+\operatorname{div}(p U)=0.
\endcases\tag1.3
$$
It follows from the transformation $p=(\g-1)P$, $t=\ds\f{T}{\g-1}$
and (1.3) that
$$
\cases
&\p_T U+\na P=0,\\
&\p_T P+P\operatorname{div} U=0.
\endcases\tag1.4
$$
Let us consider the following Cauchy problem for (1.4):
$$
\cases
&\p_T U+\na P=0,\\
&\p_T P+P\operatorname{div} U=0,\\
&U|_{T=0}=\ve U_0(x), \quad P|_{T=0}=1+\ve P_0(x),
\endcases\tag1.5
$$
where $\ve>0$ is small and $U_0(x)=(u_1^0(x),u_2^0(x))\in
C_0^{\infty}(\Bbb R^2)$, $P_0(x)\in C_0^{\infty}(\Bbb R^2)$ are
supported in the disc $B(0, M)$. One then obtains that $P$ satisfies
the nonlinear wave equation
$$\p_T\left(\f{\p_TP}{P}\right)-\Delta P=0.\tag1.6$$
Let $v(T,x)=\ln P$. Then it follows from (1.6) and the initial data in
(1.5) that
$$
\cases
&\p_T^2v-\operatorname{div}(e^v\na v)=0,\\
&v(0,x)=\ln(1+\ve P_0(x))=\ve P_0(x)+\ve^2
\dsize\sum_{n=2}^{\infty}(-1)^{n-1}\f{\ve^{n-2}}{n!}P^n_0(x),\\
&\p_tv(0,x)=-\ve \operatorname{div}U_0(x).
\endcases\tag1.7$$
In (1.7), use $t$ and $u(t,x)$ in place of $T$ and $v(T,x)$,
respectively. As a nonlinear problem equivalent to (1.7), one can then
consider
$$\cases
&\p_t^2 u-\operatorname{div}(e^u\nabla u)=0, \quad
\text{$(t,x)\in [0,\infty) \times\Bbb R^2$},\\
&u(0,x)=\varepsilon u_0(x),\\
&\partial_t u(0,x)=\varepsilon u_1(x),\\
\endcases\tag1.8
$$
where $u_0(x)=P_0(x)$ and $u_1(x)=-\operatorname{div}U_0(x)$.
In this way we have given a brief derivation on the nonlinear wave
equation in the form (1.1) from the fundamental equations of
compressible fluid dynamics.

Without loss of generality, we will assume that $c_i(0)=1$ ($i=1,2$)
in (1.1). Since third-order terms like $O(u^2D^2u)$ and $O(u|Du|^2)$
will not have an essential influence on the blowup behavior of small
data solution to problem (1.1), Eq.~(1.1) is basically equivalent to
$$\cases
&\p_t^2 u-\ds\sum_{i=1}^2\p_i\left((1+c_iu)\p_i u\right)=0,\quad
\text{$(t,x)\in [0,\infty) \times\Bbb R^2$},\\
&(u(0,x), \p_t u(0,x))=(\ve u_0(x), \ve u_1(x)),\\
\endcases\tag1.9
$$
where $c_1=2c'_1(0), c_2=2c'_2(0)$, and $c_1^2+c_2^2\not=0$.

We introduce polar coordinates $(r, \th)$ in $\Bbb R^2$,
$$
\cases
x_1=r\cos\th,\\
x_2=r\sin\th_,
\endcases
$$
where $r=\sqrt{x_1^2+x_2^2}$, $\th\in[0,2\pi]$, and $\o\equiv(\o_1,
\o_2)=(\cos\th, \sin\th)$.
Later we will need the function
$$
F_0(\si,\th)\equiv\ds\f{1}{2^{3/2}\pi}\int_{\si}^{+\infty}\ds\f{R(s,\o;
u_1)-\p_sR(s,\o; u_0)}{\sqrt{s-\si}}\,ds,\tag1.10
$$
where $\si\in \Bbb R$, and $R(s,\o;v)$ is the Radon transform of the
smooth function $v(x)$, i.e., $R(s,\o;v)=\int_{x\cdot\o=s}v(x)\,dS$.
>From Theorem 6.2.2 and (6.2.12) of [14], one has that the
function $F_0(\si,\th)\not\equiv 0$ unless $u_0(x)\equiv 0$ and
$u_1(x)\equiv 0$. Moreover, $F_0(\si,\th)\equiv 0$ for $\si\ge M$
and $\ds\lim_{\si\to -\infty}F_0(\si,\th)=0$. Therefore,
$$
  \ds\min_{\si,\th}[\p_{\si}F_0(\si,\th)(c_1\cos^2\th+c_2\sin^2\th)]<0
$$
exists as long as $(u_0(x), u_1(x))\not\equiv 0$.

We will assume throughout this paper that there is a unique point
$(\si^0, \th^0)$ such that
$$\cases &
\p_{\si}F_0(\si^0,\th^0)(c_1\cos^2\th^0+c_2\sin^2\th^0)=\ds\min_{\si\in\Bbb
R,\, \th\in [0, 2\pi]}
[\p_{\si}F_0(\si,\th)(c_1\cos^2\th+c_2\sin^2\th)],\\
&\text{the Hessian matrix }\na_{\si, \th}^2[\p_{\si}F_0(\si,
\th)(c_1\cos^2\th+c_2\sin^2\th)]|_{(\si, \th)=(\si^0, \th^0)}>0.
\endcases\tag1.11
$$

Let $T_{\ve}$ denote the lifespan of the smooth solution to
(1.9). Then one has:

\proclaim{Theorem 1.1} Let $u_0(x)$, $u_1(x)\in C_0^{\infty}(\Bbb R^2)$
%%and $(u_0(x), u_1(x))\not\equiv 0$ which are
be supported in the disc $B(0, M)$ and let assumption\/
\rom{(1.11)} hold. Then\rom:
\roster
\item
$$
\ds\lim_{\ve\rightarrow 0}\ve\sqrt{T_\ve}=\tau _0\equiv-\f{1}{\p_{\si}
F_0(\si^0,\th^0)(c_1\cos^2\th^0+c_2\sin^2\th^0)}>0.\tag 1.12
$$

\item There exists a point $M_\ve=(T_\ve, x_\ve)$ and a positive
constant $C$ independent of $\ve$ such that
\widestnumber\item{(iiiiiii)}
\item"(i)" $u(t,x)\in C([0,T_\ve]\times\Bbb R^2)$ and
$\|u\|_{L^{\infty}([0,T_\ve]\times\Bbb R^2)}\leq C\ve$\rom.

\item"(ii)" $u\in C^2([0,T_\ve]\times\Bbb R^2\setminus\{M_{\ve}\})$,
and, for $t <T_\ve$, it satisfies
$$
\align
\|\nabla_{t, x} u(t,\cdot)\|_{L^\infty(\Bbb
R^2)}&\leq\f{C}{T_\ve-t},\tag1.13 \\
\|\p_t u(t,\cdot)\|_{L^\infty(\Bbb
R^2)}&\geq\f{1}{C(T_\ve-t)}.\tag1.14
\endalign
$$
\endroster
\endproclaim

\remark{Remark \rom{1.1}} Compared with the ``lifespan theorems'' of
[4-5], Theorem~1.1 states that the solution $u(t,x)$ to (1.9) is
continuous up to the blowup time $t=T_{\ve}$, while its first-order
derivatives $\nabla_{t, x} u$ develop a singularity at $t=T_{\ve}$. In
the terminology of [4-5], this corresponds to an ``ODE blowup.'' On
the contrary, the blowup result of [4-5] on small data solutions to
the 2-D nonlinear wave equation $\p_t^2v-\Delta_xv +\ds\sum_{0\le
i,j,k\le 2}g_{ij}^k\p_k v\,\p_{ij}^2v=0$, where the nonlinearity
depends on the derivatives of $v$, but not $v$ itself, shows that the
solution $v(t,x)$ is $C^1$ up to the blowup time $T_{\ve}$, while the
second-order derivatives $\na_{t,x}^2v$ develop a singularity at
$t=T_{\ve}$. In the terminology of [4-5], this is a ``geometric
blowup.'' \endremark

\remark{Remark \rom{1.2}} One readily obtains $u(t,x)\in
C^{\infty}(([0, T_\ve]\times\Bbb R^2)\setminus\{M_{\ve}\})$ from
$u(t,x)\in C^2(([0, T_\ve]\times\Bbb R^2)\setminus\{M_{\ve}\})$ in
Theorem 1.1. Since $u(t,x)\in C^{\infty}([0, T_{\ve})\times\Bbb R^2)$,
this in fact follows from the property of finite propagation speed
which holds for hyperbolic equations. \endremark

\remark{Remark \rom{1.3}} In particular, for the 2-D pressure-gradient
model $\p_t^2 u -\ds\sum_{i=1}^2\p_i(e^u\p_i u)=0$ with small initial
data $(u(0,x), \p_tu(0,x))=(\ve u_0(x), \ve u_1(x))$ and $(u_0(x),
u_1(x))\not\equiv 0$, it follows from Theorem 1.1 that the lifespan
$T_{\ve}$ of the smooth solution $u(t,x)$ satisfies $\ds\lim_{\ve\to
0}\ve\,\sqrt{T_{\ve}} =-\ds\f{1}{\ds\min_{\si\in\Bbb R,\, \th\in [0,
2\pi]}\p_{\si}F_0(\si,\th)}$ under an assumption on the function
$\p_{\si}F_0(\si,\th)$ that is analogous to (1.11). We thus have
extended the blowup result of [21] valid for the rotationally
symmetric case to this now more general situation. In addition,
returning to the original pressure-gradient system (1.5), one obtains
that $\p_tP$ and $\operatorname{div} U$ develop a singularity at time
$t=T_{\ve}$. This corresponds to the formation of a shock emanating
from the blowup point as shown in [26] for the compressible Euler
system. \endremark

\remark{Remark \rom{1.4}} The nonlinear equation (1.9) can be
rewritten as $\p_t^2 u-(1+u)\Delta u=|\nabla u|^2$ when
$c_1=c_2=1$. For the 3-D equation $\p_t^2 u-(1+u)\Delta u=0$ with
small initial data $(u(0,x), \p_t u(0,x))=(\ve u_0(x), \ve u_1(x))$,
in [6, 23] it was shown that smooth solutions exist globally. On the
other hand, for the $n$-dimensional nonlinear wave equation ($n=2,3$)
with coefficients depending on the derivatives of the solution,
$\p_t^2 u-c^2(\p_t u)\Delta u=0$ and, more generally,
$\dsize\sum_{i,j=0}^ng_{ij}(\na u)\p_{ij}^2u=H(\na u)$, where $t=x_0$,
$x=(x_1, ..., x_n)$, $g_{ij}(\na u)=c_{ij}+O(|\na u|)$, $H(\na
u)=O(|\na u|^2)$, and the linear part
$\dsize\sum_{i,j=0}^nc_{ij}\p_{ij}^2u$ is strictly hyperbolic with
respect to time $t$, it is known that small data smooth solutions
exist globally if related null conditions hold (see [8, 14] and
others), while otherwise small data smooth solutions blow up in finite
time (see [4-5, 13, 17, 22] and others). We point out that in the case
considered here the coefficients of the nonlinear equation (1.9)
depend on both the solution $u$ and its derivatives.
\endremark

\medskip

Near the blowup point $M_{\ve}$ one can give a more accurate
description of the behavior of the solution $u(t,x)$ which is similar
to statements in the ``geometric blowup theorems'' of [4--5].

\proclaim{Theorem 1.2} Assume that the constants $\tau_1$, $A_0$,
$A_1$ and $\dl_0$ satisfy $0<\tau_1<\tau_0$, $A_0<\si^0<A_1<M$ and
that $\dl_0>0$ is sufficiently small. Moreover, assume that $A_0$ and
$A_1$ are close to $\si^0$. Denote by $D$ the domain
$$
D\equiv\{(s,\th,\tau)\mid A_0\le s\le A_1, \,
\th^0-\dl_0\le\th\le\th^0+\dl_0, \, \tau_1\le\tau\le\tau_\ve\},
$$
where $\tau_\ve=\ve\sqrt{T_\ve}$. Then there exist a subdomain $D_0$
of $D$ containing a point $m_\ve=(s_\ve,\th_\ve,\tau_\ve)$ and
functions $\phi(s,\th,\tau), \,v(s,\th,\tau)\in C^3(D_0)$ with the
following properties\rom:
\roster
\item In the domain $D_0$, $\phi$ satisfies
$$
\cases
\enspace \p_s\phi(s,\th,\tau)\geq 0,\quad\p_s\phi(s,\th,\tau)=
0\Longleftrightarrow(s,\th,\tau)=m_\ve,\\
\enspace \p_{\tau s}^2\phi(m_\ve)<0,\quad\nabla_{s,\th}\p_s\phi(m_\ve)=0,
\quad\nabla^2_{s,\th}\p_s\phi(m_\ve)>0.
\endcases\tag H
$$

\item It holds
$$
v(m_\ve)\neq 0.\tag1.15
$$
\endroster
Moreover, let the function $G(\si,\th,\tau)$ be defined by
$G(\Phi)=v(s,\th,\tau)$ in the domain $\Phi(D_0)$, where $\Phi$ is a
map such that $\Phi(s,\th,\tau)=(\phi(s,\th,\tau),\th,\tau)$. Then
$u(t,x)=\ds\f{\ve}{\sqrt{r}}\,G(r-t,\th,\ve\sqrt{t})$
solves \rom{Eq.~(1.9)}.
\endproclaim

\remark{Remark \rom{1.5}} Theorem~1.2 provides a more accurate
description of the solution near the blowup point $M_{\ve}
=\Phi(m_{\ve})$ than Theorem~1.1. First, one has that
$G(\si,\th,\tau)\in C(\Phi(D_0))$ because of $\phi,v\in C^3(D_0)$ and
(H) of Theorem~1.2. To prove this assertion, we are only required to
show that $G$ is continuous at the point $M_{\ve}=(\si_\ve,
\th_\ve, \tau_{\ve})\equiv (\phi(m_\ve), \th_\ve, \tau_{\ve})$. To
this end, let $(\si_n, \th_n, \tau_n)\in
\Phi(D_0)$ be such that $(\si_n, \th_n, \tau_n)\to
(\si_\ve,\th_\ve,\tau_\ve)$ as $n\to\infty$. It then follows
from (H) that there is a unique point $(s_n, \th_n,
\tau_n)\in D_0$ such that $\si_n=\phi(s_n, \th_n, \tau_n)$. By
Taylor's formula, one has $\si_n-\si_\ve=\na_{\th,\tau} \phi(m_{\ve})\cdot
(\th_n-\th_\ve, \tau_n-\tau_\ve)+\ds\f12\,\p_{s\tau}^2\phi(m_{\ve})(s_n-s_\ve)
(\tau_n-\tau_\ve)+\ds\f12\,(\th_n-\th_\ve,
\tau_n-\tau_\ve)\na^2_{\th,\tau}\phi(m_\ve)(\th_n-\th_\ve,
\tau_n-\tau_\ve)^T+\ds\f16\,\p_s^3\phi(m_\ve)(s_n-s_\ve)^3
+o(|s_n-s_\ve|^3)+o(|\th_n-\th_\ve| +|\tau_n-\tau_\ve|)$. Together
with $\p_s^3\phi(m_\ve)>0$, this yields $s_n\to s_\ve$ as
$n\to\infty$. Therefore, one obtains $G\in C(\Phi(D_0))$ from
$G(\Phi)=v$ and the continuity of $v, \phi$ in $D_0$. It follows that
$u(t,x)=\ds\f{\ve}{\sqrt{r}}\,G(r-t,\th,
\ve\sqrt t)\in C\left(\left[\f{\tau_1^2}{\ve^2}, T_{\ve}\right]\times
\Bbb R^2\right)\cap
C^1\left(\left(\left[\f{\tau_1^2}{\ve^2}, T_{\ve}\right]\times\Bbb
R^2\right)\setminus\{M_{\ve}\}\right)$ and $\|u\|_{L^{\infty}}\le
C\ve^2$ in a neighborhood of $M_{\ve}$. Regarding the other properties
of $u(t,x)$ near $M_{\ve}$ stated in Theorem~1.1, see \S4 below for
details.
\endremark

\medskip

There are some interesting papers on the Riemann problem for the
pressure-gradient system (1.5) and (1.6), respectively, with special
discontinuous initial data, with either a mathematical treatment or a
numerical simulation (see [1, 19-20, 25, 29-31] and the references
therein). There are also many results on the blowup of classical
solutions and the global existence and uniqueness of weak solutions,
respectively, to 1-D variational wave equations (see [2, 7, 10-12, 16,
27-28] and the references therein). In the multidimensional case of
Eq.~(1.1), however, except for the rotationally symmetric case, where
in [9, 21] blowup results have been established, until now there were
no results on the finite-time blowup of smooth solutions to (1.1) or
even on mechanisms of this blowup. In this paper, we shall focus on
these two problems, i.e., we will establish the precise lifespan
$T_{\ve}$ in Theorem~1.1 and determine the blowup mechanism in
Theorem~1.2.

Let us comment on the proofs of Theorems~1.1 and~1.2. First we derive
the required lower bound on the lifespan $T_{\ve}$ for solutions to
problem (1.9).  As in [14, Chapter~6] and [13], by constructing a
suitable approximate solution $u_a(t,x)$ to (1.9) and then considering
the difference of the exact solution $u(t,x)$ and $u_a(t,x)$, applying
the Klainerman-Sobolev inequality, and further establishing a delicate
energy estimate, we obtain this lower bound on the lifespan
$T_{\ve}$. Next we derive the required upper bound on
$T_{\ve}$. Motivated by the ``geometric blowup'' method of [4-5], we
introduce the blowup system of (1.9) to study simultaneously the
lifespan $T_{\ve}$ and blowup mechanism of smooth solution $u$. That
is, by introducing a singular change of coordinates $\Phi$ in the
domain $D=\left\{(\si, \th, \tau)\mid -C_0\le\si\le M,\, 0\le\th\le
2\pi,\, 0<\tau_1\le\tau\le \tau_{\ve} \right\}$,
$$
\text{$(s, \th, \tau)\to (\phi(s,\th,\tau),\th, \tau)$, where
$\phi(s,\th,\tau_1)=s$ and
$\p_s\phi=0$ holds at some point,}
$$
where $\si=r-t$, $\tau=\ve\sqrt{t}$, and $C_0>0$ a fixed constant, and
setting $G(\Phi)=v(s,\th,\tau)$, we obtain a nonlinear system for
$(\phi, v)$ from the ansatz $u(t,x)=\ds\f{\ve}{\sqrt{r}}\,
G(r-t, \th, \ve \sqrt{t})$ and the equation in (1.9). This blowup
system for (1.9) has a unique smooth solution $(\phi, v)$ for
$\tau\le\tau_{\ve}$, where the couple $(\phi, v)$ satisfies properties
(H) and (1.15) of Theorem~1.2. This enables us to determine the blowup
point at time $t=T_{\ve}$ for the solution $u$ of (1.9) and give a
complete asymptotic expansion of $T_{\ve}$ as well as a precise
description of the behavior of $u(t,x)$ close to the blowup point. In
order to treat the resulting blowup system, as in [4-5], we use the
Nash-Moser-H\"ormander iteration method to overcome the difficulties
introduced by the free boundary $t=T_{\ve}$ and the inherent
complexity of the nonlinear blowup system. To this end, the linearized
system is solved first. Thanks to the energy estimates established in
[4-5], we are then able to complete the proof of Theorem~1.2.

\smallskip

The paper is organized as follows: In \S2, as in [9, 21], we construct
a suitable approximate solution $u_a(t,x)$ to (1.9) and establish
related estimates, which allows us to obtain the required lower bound
on the lifespan $T_{\ve}$. In \S3, the blowup system for (1.9) is
solved, which allows us to prove Theorem~1.2. Then, in $\S 4$, we
conclude the proof of Theorem~1.1 based on Theorem 1.2.

\subsubhead Notation \endsubsubhead
Throughout the paper, we will use the following notation:
$Z$ denotes one of the Klainerman vector fields in $\Bbb
R_t^+\times\Bbb R^2$, i.e.,
$$
\p_t, \enspace \p_i, \enspace \G_0=t\p_t+\ds\sum_{j=1}^2x_j\p_j,
\enspace H_i=x_i\p_t+t\p_i, \enspace i=1, 2, \enspace  R=x_1\p_2-x_2\p_1,
$$
$\p$ stands for $\p_t$ or $\p_i$ ($i=1,2$), and $\na_x$ stands for
$(\p_1, \p_2)$.

%% ------------------------------------------------------------------------

\head \S2. Lower bound on the lifespan $T_\varepsilon$\endhead

In this section, we establish the lower bound of
$T_\varepsilon$ for smooth solution to the Cauchy problem (1.9).

Let $\tau=\varepsilon\sqrt{1+t}$ be the slow time variable and assume
the solution to (1.9) can be approximated by
$$
\f{\varepsilon}{\sqrt{r}}\,V(\si,\th, \tau),\qquad r>0,
$$
where $\si=r-t$, $(x_1, x_2)=(r\cos\th, r\sin\th)$ with $\th\in [0,
2\pi]$.

The function $V(\si, \th, \tau)$ solves the equation
$$
\cases
&\p_{\si\tau}^2 V+(c_1\cos^2\th+c_2\sin^2\th)V\p_\si^2 V+
(c_1\cos^2\th+c_2\sin^2\th)(\p_\si V)^2=0,\\
V(\si,\th, 0)=F_0(\si,\th) ,\\
&\operatorname{supp} V(\cdot, \th, \tau)\subseteq \{\si\leq M\},\\
\endcases\tag2.1
$$
where $F_0(\si, \th)$ has been defined in (1.10).

For problem (2.1), one has:

\proclaim{Lemma 2.1}  \rom{Eq.~(2.1)} admits a $C^\infty$ solution for
$0\leq\tau<\tau_0$ with the number $\tau_0$ being given
in \rom{(1.12)}.
\endproclaim

\demo{Proof}  Set $U(\si,\th,\tau)=\p_\si V(\si,\th,\tau)$. Then it
follows from (2.1) that
$$\cases
&\p_\tau U+(c_1\cos^2\th+c_2\sin^2\th)V
\p_\si U+(c_1\cos^2\th+c_2\sin^2\th)U^2=0,\\
&U(\si,\th,0)=\p_\si F_0(\si,\th).\\
\endcases\tag 2.2
$$

The characteristic curve $\si=\si(s,\th, \tau)$ of (2.2) starting at
the point $(s, \th, 0)$ is defined by
$$\cases
&\ds\f {d\si}{d\tau}(s,\th, \tau)=(c_1\cos^2\th+c_2\sin^2\th)
V(\si(s,\th, \tau),\th,\tau),\\
&\si(s,\th, 0)=s.\\
\endcases\tag2.3
$$
Along this characteristic curve, it follows from (2.2) that, for
$\tau<\tau_0$,
$$
U(\si(s,\th, \tau),\th,\tau)=\ds\f{\p_\si F_0(s,\th)}{1+
(c_1\cos^2\th+c_2\sin^2\th)\p_\si F_0(s,\th)\tau}.\tag2.4
$$
Because of $U(\si,\th,\tau)=\p_\si V(\si(s,\th, \tau),\th, \tau)$,
from (2.3)-(2.4) one then obtains that
$$
\cases
&\p_{\tau s}^2\si(s,\th,
\tau)=\ds\f{(c_1\cos^2\th+c_2\sin^2\th)\p_\si
F_0(s,\th)}{1+(c_1\cos^2\th
+c_2\sin^2\th)\p_\si F_0(s,\th)\tau}\,\p_s\si(s,\th, \tau),\\
&\p_s\si(s,\th, 0)=1.\\
\endcases
$$
This gives $\partial_s\si(s,\th,
\tau)=1+(c_1\cos^2\th+c_2\sin^2\th)\p_\si F_0(s,\th)\tau>0$ for
$0\leq \tau<\tau_0$ and then
$$
\si(s,\th, \tau)=\si(M,\th, \tau)+s-M+(c_1\cos^2\th+c_2\sin^2\th)
F_0(s,\th)\tau\tag2.5
$$
and
$$
V(\si(s,\th, \tau),\th,\tau)=\f{\p_\tau \si(M,\th,
\tau)}{c_1\cos^2\th+c_2\sin^2\th}+F_0(s,\th).\tag2.6
$$
Note that $\si(M,\th, \tau)=M$ such that $V(\si,\th,\tau)$ satisfies
the boundary condition $V|_{\si=M}=0$. This, together with
(2.5)-(2.6), yields $V(\si,\th,\tau)=F_0(s,\th)$ and
$\si=s+(c_1\cos^2\th+c_2\sin^2\th)F_0(s,\th)\tau$. By the implicit
function theorem, one then has that $s=s(\si,\th,\tau)$ is a smooth
function of $\si,\th,\tau$ for $\tau<\tau_0$. Therefore,
$V(\si,\th,\tau)=F_0(s(\si,\th,\tau),\th)$ is a smooth solution of
(2.1) for $0\leq \tau<\tau_0$ as claimed. \qed
\enddemo

>From [14, Chapter~6], one has that $F_0(\si,\th)\in C^{\infty}(\Bbb
R)$ is supported in $(-\infty,M]$ and obeys the estimates
$$
\left|\p_{\si}^k\p_{\th}^lF_0(\si,\th)\right|\leq
C_{kl}\,(1+|\si|)^{-1/2-k},\quad
k\in \Bbb N_0.\tag 2.7
$$

>From (2.7), we now derive a decay estimate of $V(\si,\th,\tau)$ in
(2.1) for $\tau<\tau_0$ and $\si\to -\infty$.

\proclaim{Lemma 2.2} For any positive constant $b<\tau_0$, one has that, in the
domain
$$
  \{(\si, \th, \tau)\mid -\infty<\si\le M, 0\le\th\le 2\pi,
  0\leq\tau\leq b\},
$$
and for $r\geq t/3$, the smooth solution $V$
to \rom{(2.1)} obeys the estimates
$$
|Z^{\al}\p_{\tau}^{l}\p_{\si}^{m}V(\si,\th,\tau)|\leq C_{\al
b}^{lm}(1+|\si|)^{-1/2-l-m},\quad \al,l,m\in \Bbb N_0,\tag2.8
$$
where $C_{\al b}^{lm}$ are positive constants depending on $b$ and
$\al, l, m$.
\endproclaim

\demo{Proof} When $\tau\leq b$, it follows from (2.5) and the
support property of $F_0(\si, \th)$ that $\ds\f{|s|}{2}\leq
|\si|\leq 2|s|$ for large $|s|$. Together with (2.6), this yields
$$
|V(\si,\th,\tau)|\leq C_b\,(1+|\si|)^{-1/2},\quad
|\p_{\si}V(\si,\th,\tau)|\leq C_b\,(1+|\si|)^{-3/2}.\tag 2.9
$$

By (2.6) and (2.4), one has
$$
\p_{\si}s(\si,\th,\tau)=\f{1}{1+\left(c_1\cos^2\th+c_2\sin^2\th\right)
\p_sF_0(s,\th)\tau}
$$
and
$$\quad\quad\p_{\si}^2
 V(\si(s,\th, \tau),\th,\tau)
=\f{\p_s^2F_0(s,\th)}{(1+\left(c_1\cos^2\th+c_2\sin^2\th\right)
\p_sF_0(s,\th)\tau)^3},
$$
which yields
$$|\p_{\si}^2 V(\si,\th,\tau)|\leq C_b\,(1+|\si|)^{-5/2}.\tag
2.10$$

Further, it follows from (2.1) and (2.10) that
$$
|\p_{\tau\si}^2V(\si,\th,\tau)|\leq C_b\,(1+|\si|)^{-3}
$$
and then
$$
|\p_{\tau}V(\si,\th,\tau)|\leq C_b\,(1+|\si|)^{-2}.\tag 2.11
$$

Based on (2.9)--(2.11), by an inductive argument one arrives at
$$
|\p_{\tau}^{l}\p_{\si}^{m}V(\si,\th,\tau)|\leq
C_b^{lm}\,(1+|\si|)^{-1/2-l-m},\quad l,m\in\Bbb N_0.
$$

Because of
$$
\align
\G_0&=\si\p_{\si}+\ds\f{\ve t}{2\sqrt{1+t}}\,\p_{\tau}, \quad
H_1=-\si\cos\th\p_{\si}+\ds\f{\ve
x_1}{2\sqrt{1+t}}\,\p_{\tau}-\f{x_2t}{r^2}\,\p_\th, \\
H_2&=-\si\sin\th\p_{\si}+\ds\f{\ve
x_2}{2\sqrt{1+t}}\,\p_{\tau}+\f{x_1t}{r^2}\,\p_\th, \quad
R=\p_\th,
\endalign
$$
one analogously obtains
$$
|Z^{\al}\p_{\tau}^{l}\p_{\si}^{m}V(\si,\th,\tau)|\leq C_{\al b}^{lm}\,
(1+|\si|)^{-1/2-l-m},\quad \al,l,m\in\Bbb N_0,
$$
which completes the proof of Lemma~2.2.
\qed
\enddemo

Next, we construct an approximate solution $u_a$ to (1.9) for
$0\leq\tau=\ve\sqrt{1+t}<\tau_0$.

Let $w_0$ be the solution of the linear wave equation
$$
\cases
&\p_t^2 w_0-\triangle w_0=0,\\
&w_0(0,x)=u_0(x),\quad \p_t w_0(0,x)=u_1(x).
\endcases
$$

It follows from [14, Theorem~6.2.1] that, for any constants $l>0$
and $0<m<1$,
$$\align
&\left|Z^{\al}(w_0(t,x)-r^{-1/2}F_0(\si,\th))\right|\leq C_{\al
l}\,(1+t)^{-3/2}(1+|\si|)^{1/2},\quad r\geq lt, \tag 2.12\\
&\left|\p^{k}w_0(t,x)\right|\leq C_{km}\,(1+t)^{-1-|k|},\quad r\le mt.\tag 2.13
\endalign$$
Choose a $C^\infty$ function $\chi(s)$ such that $\chi(s)=1$ for
$s\leq 1$ and $\chi(s)=0$ for $s\geq 2$. For $0\leq
\tau=\ve\sqrt{1+t}<\tau_0$, we take the approximate solution $u_a$
to (1.9) to be
$$
u_a(t,x)=\varepsilon\left(\chi(\varepsilon
t)w_0(t,x)+r^{-1/2}(1-\chi(\varepsilon
t))\chi(-3\ve\si)V(\si,\tau)\right).\tag 2.14
$$
By Lemma 2.2 and [14, Theorem 6.2.1], one has that, for a fixed
positive constant $b<\tau_0$,
$$|Z^\al
u_a(t,x)|\leq C_{\al b}\,\ve
(1+t)^{-1/2}(1+|\si|)^{-1/2},\quad \tau\leq b.\tag 2.15
$$

Set $J_a=\p_t^2 u_a-(1+c_1u_a)\p_1^2 u_a-(1+c_2u_a)\p_2^2
u_a-c_1(\p_1u_a)^2-c_2(\p_2u_a)^2$.

\proclaim{Lemma 2.3} One has
$$\ds \int_0^{b^2\!/\ve^2-1}\|Z^\alpha
J_a(t,\cdot)\|_{L^2}\,dt\leq C_{\al b}\,\ve^{3/2}.
$$
\endproclaim

\demo{Proof} We divide the proof into three parts.
\enddemo

\subhead (A) $0\leq t\leq \ds\f{1}{\varepsilon}$ \endsubhead
In this case, $\chi(\ve t)=1$ and $u_a=\ve w_0$.  This yields
$$
J_a=-\ve^2w_0c_1\p_1^2 w_0-\ve^2w_0c_2\p_2^2 w_0
-\ve^2c_1(\p_1w_0)^2-\ve^2c_2(\p_2w_0)^2.
$$
It follows from (2.15) and a direct computation that, for $0\leq
t\leq \ds\f{1}{\ve}$,
$$
\|Z^{\al}J_a(t,\cdot)\|_{L^2}\leq C\,\ve^2 (1+t)^{-1/2}.\tag 2.16
$$

\subhead (B) $\ds\f{1}{\varepsilon}\leq t\leq\ds\f{2}{\varepsilon}$
\endsubhead
We now rewrite $u_a$ as
$$
u_a=\ve w_0(t,x)+\ve(1-\chi(\ve
t))\left(r^{-1/2}\chi(-3\ve\si)V(\si,\th,\tau)-w_0(t,x)\right).
$$
Then
$$J_a=J_1+J_2+J_3+J_4,\tag 2.17$$
where
$$
\align
J_1&=-c_1u_a\p_1^2u_a-c_2u_a\p_2^2u_a-c_1
(\p_1u_a)^2-c_2(\p_2u_a)^2,\\
J_2&=\ve(\p_t^2-\Delta)\left[(1-\chi(\ve
t))r^{-1/2}\chi(-3\ve\si)\bigl(V(\si,\th,\tau)-F_0(\si,\th)\bigr)
\right],\\
J_3&=\ve(\p_t^2-\Delta)\left[(1-\chi(\ve
t))\chi(-3\ve\si)
\bigl(r^{-1/2}F_0(\si,\th)-w_0(t,x)
\bigl)\right],\\
J_4&=\ve(\p_t^2-\Delta)\left[(1-\chi(\ve
t))(\chi(-3\ve\si)-1)w_0(t,x)\right].
\endalign$$

We treat each term $J_i$ ($1\le i\le 4$) in (2.17) separately.
>From (2.15) one obtains
$$\|Z^{\al}J_1(t, \cdot)\|_{L^2}\leq C_{\al b}\,\ve^2(1+t)^{-1/2}.\tag 2.18$$
Since
$$
\align
J_2&=\ve r^{-1/2}(\p_t-\p_r)(\p_t+\p_r)\left[(1-\chi(\ve
t))\chi(-3\ve\si)\bigl(V(\si,\th,\tau)-F_0(\si,\th)\bigr)\right]\\
&\qquad-\f{\ve}{4}\,r^{-5/2}(1-\chi(\ve
t))\chi(-3\ve\si)\bigl(V(\si,\th,\tau)-F_0(\si,\th)\bigr)\\
&\qquad-\ve r^{-5/2}\p_\th^2\left[(1-\chi(\ve
t))\chi(-3\ve\si)\bigl(V(\si,\th,\tau)-F_0(\si,\th)\bigr)\right]
\endalign
$$
and $V(\si,\th,\tau)-F_0(\si,\th)=\int_0^\tau \p_\tau V(\si,\th,s)ds$, one has
$$
\|Z^{\al}J_2(t,\cdot)\|_{L^2}\leq C_{\al
b}\,\ve^2(1+t)^{-1/2}.\tag 2.19
$$
Note that $-\ds\f{2}{3\ve}\leq\si\leq M$ holds on the support of $J_3$
which implies $r\geq \ds\f{1}{3}t$. This, together with (2.12), yields
$$
\|Z^{\al}J_{3}(t, \cdot)\|_{L^2}\leq
C_{\al}\,\ve^2(1+t)^{-1/2}.\tag2.20
$$
Analogously, together with (2.13), one arrives at
$$
\|Z^{\al}J_4(t, \cdot)\|_{L^2}\leq C_{\al b}\,\ve^2(1+t)^{-1}.\tag 2.21
$$
Collecting (2.18)-(2.21) yields
$$
\|Z^{\al}J_a(t, \cdot)\|_{L^2}\leq C_{\al b}\,\ve^2(1+t)^{-1/2},\quad
\f{1}{\ve}\leq t\leq
\f{2}{\ve}.\tag 2.22
$$

\subhead (C) $\ds\f{2}{\varepsilon}\leq t\leq \f{b^2}{\ve^2}-1$ \endsubhead
A direct computation yields
$$
\align
J_a&=-\ve^2 r^{-1/2}\p_{\tau\si}^2\hat
V\biggl(\f{1}{\sqrt{1+t}}-r^{-1/2}\biggr)\\
&\qquad-\ve^2
r^{-1}\biggl(\p_{\tau\si}^2\hat V+(c_1\cos^2\th+c_2\sin^2\th)\hat V\p_{\si}^2\hat
V+(c_1\cos^2\th+c_2\sin^2\th)(\p_{\si}\hat V)^2\biggr)\\
&\qquad+O(\ve^2)\,(1+t)^{-3/2}(1+|\si|)^{-1/2},\tag 2.23
\endalign
$$
where $\hat V(\si,\th,\tau)=\chi(-3\ve\si)V(\si,\th,\tau)$.
It follows from (2.1) that
$$
\multline
\ve^2 r^{-1}  \biggl(\p_{\tau\si}^2
\hat V+(c_1\cos^2\th+c_2\sin^2\th)\hat V\p_{\si}^2\hat V
+(c_1\cos^2\th+c_2\sin^2\th)(\p_{\si}\hat V)^2\biggr)\\
= O(\ve^2)\,(1+t)^{-3/2}(1+|\si|)^{-3/2},
\endmultline
\tag 2.24
$$
here we have used the fact that $\chi(-3\ve\si)(1-\chi(-3\ve\si)$ is
supported in the interval $[-\ds\f{2}{3\ve}, -\ds\f{1}{3\ve}]$.
Substituting (2.24) into (2.23) yields
$$\|Z^{\al}J_a(t, \cdot)\|_{L^{2}}\leq C_{\al
b}\ve^2(1+t)^{-3/4}.\tag 2.25
$$

\medskip

Consequently, combining (2.16), (2.22), and (2.25) yields
$$
\ds\int_0^{b^2\!/\ve^2-1}\|Z^\al J_a(t, \cdot)\|_{L^2}\,dt\leq
C_{\al b}\,\ve^{3/2},
$$
which completes the proof of Lemma 2.3.
\qed

\proclaim{Lemma 2.4} For sufficiently small $\varepsilon$ and
$0\leq \tau=\ve\sqrt{1+t}\leq b<\tau_0$, \rom{Eq.~(1.9)} admits a
$C^\infty$ solution $u$ which satisfies the estimate
$$
|Z^\kappa\p (u-u_a)|\leq C_{b}\,
\varepsilon^{3/2}(1+t)^{-1/2}(1+|t-r|)^{-1/2}.\tag2.26
$$
for $|\kappa|\le 2$.
\endproclaim

\demo{Proof} Set $v=u-u_a$. Then
$$\cases
&\p_t^2 v-(1+c_1u)\p_1^2v-(1+c_2u)\p_2^2v=F,\\
&v(0,x)=\p_t v(0,x)=0,
\endcases\tag 2.27
$$
where
$$
F=-J_a+c_1v\p_1^2u_a+c_2v\p_2^2u_a
+c_1(\p_1v)^2+c_2(\p_2v)^2
+2c_1(\p_1 v)(\p_1u_a)
+2c_2(\p_2 v)(\p_2u_a).\tag 2.28
$$
We will use continuous induction to prove (2.26). To this end, we
assume that, for some $T\leq \ds\f{b^2}{\ve^2}-1$,
$$
|Z^{\kappa}\p
v|\leq\varepsilon\,(1+t)^{-1/2}(1+|t-r|)^{-1/2},\quad
|\kappa|\leq 2,\enspace t\leq T, \tag 2.29
$$
holds and subsequently we prove that
$$
|Z^{\kappa}\p
v|\leq\f{\varepsilon}{2}\,(1+t)^{-1/2}(1+|t-r|)^{-1/2},\quad
|{\kappa}|\leq 2, \enspace t\leq T.\tag2.30
$$
Note that from (2.29) one has
$$|Z^{\kappa} v|\leq
C\varepsilon\,(1+t)^{-1/2}(1+|t-r|)^{1/2},\quad |{\kappa}|\leq
2, \enspace t\leq T.\tag2.31
$$
Applying $Z^\alpha$ to both hand sides of (2.27) yields, for
$|\alpha|\leq 4$,
$$
\multline
\bigl(\p_t^2 -(1+c_1u)\p_1^2-(1+c_2u)\p_2^2\bigr)Z^\alpha v=G\\
\equiv\sum_{|\beta|\leq
|\al|}C_{\al\beta}Z^{\beta}F+\bigl[Z^{\al},c_1u\p_1^2+c_2u\p_2^2)]v\quad
+\sum_{|\beta|<|\al|}C_{\al\beta}'Z^{\beta}\bigl(c_1u\p_1^2v+c_2u\p_2^2v\bigr),
\endmultline
\tag2.32
$$
where the commutator relation $[Z^{\al}, \p_t^2-\triangle]
=\dsize\sum_{|\beta|<|\al|}C''_{\al\beta}Z^{\beta}(\p_t^2-\triangle)$
was used, and $C_{\al\beta}$, $C'_{\al\beta}$, $C''_{\al\beta}$ are
suitable constants.

Next we derive from (2.32) an estimate of $\|\p Z^\alpha
v(t, \cdot)\|_{L^2}$. Define the energy
$$
E(t)=\f{1}{2}\sum_{|\al|\leq 4}\int_{\Bbb R^2}(|\p_t Z^\alpha
v|^2+(1+c_1u)(\p_1 Z^\al v)^2+(1+c_2u)(\p_2 Z^\al v)^2)\,dx.
$$
Multiplying both sides of (2.32) by $\p_t Z^\alpha v$ ($|\al|\le 4$),
integrating by parts in $\Bbb R^2$, and noting that $|\p u|=|\p u_a+\p
v|\leq C_b\ve(1+t)^{-1/2}$ from the construction of $u_{a}$ and
assumption (2.29), one arrives at
$$
E'(t)\leq \f{C_b\ve}{\sqrt{1+t}}\,E(t)+\sum_{|\al|\leq 4}\int_{\Bbb
R^2}|G|\cdot|\p_t Z^\al v|\,dx.\tag2.33
$$
Moreover, due to the inductive hypothesis (2.29) and (2.15), one
has
$$|Z^{\kappa}u|\leq C_{b}\ve\,(1+t)^{-1/2}(1+|\si|)^{1/2}\leq C_{b}\ve,
\quad |\kappa|\leq 2, \enspace t\leq T.\tag 2.34$$

We now treat each term in the sum $\dsize\sum_{|\al|\leq 4}\int_{\Bbb
R^2}|G|\cdot|\p_t Z^\al v|\,dx$ separately.
\enddemo

\subhead (A) Estimation of $\dsize\sum_{|\beta|<|\al|} \int_{\Bbb
R^2}|Z^{\beta}\bigl(c_1u\p_1^2v+c_2u\p_2^2v\bigr)|\cdot|\p_t Z^\alpha
v|\,dx$ \endsubhead
It follows from (2.34) that, for $|\beta|<|\al|$, $i=1,2$,
$$
\align
\int_{\Bbb R^2}|Z^{\beta} (u\p_i^2 v)|\cdot|\p_t
Z^{\al}v|\, dx
&\leq C_b\sum_{|\beta_1|+|\beta_2|=|\beta|}\int_{\Bbb
R^2}|Z^{\beta_1}u|\cdot|Z^{\beta_2}\p_i^2 v|\cdot|\p_t
Z^{\al}v|\,dx\\
&\leq C_b\sum_{|\beta_1|+|\beta_2|=|\beta|}\int_{\Bbb
R^2}|Z^{\beta_1}v|\cdot|Z^{\beta_2}\p_i^2 v|\cdot|\p_t
Z^{\al}v|\,dx \tag2.35\\
&\qquad +C_b\sum_{|\beta_1|+|\beta_2|=|\beta|}\int_{\Bbb
R^2}|Z^{\beta_1}u_a|\cdot|Z^{\beta_2}\p_i^2 v|\cdot|\p_t
Z^{\al}v|\,dx.
\endalign
$$
Due to $$
\p_t=\ds\f{t\G_0-\ds\sum_{i=1}^2x_iH_i}{t^2-r^2}, \quad
\p_1=\ds\f{x_2R+tH_1-x_1\G_0}{t^2-r^2}, \quad
\p_2=\ds\f{-x_1R+tH_2-x_2\G_0}{t^2-r^2},
$$
one then has
$$
|Z^{\beta_2}\p_i^2 v|\leq
\f{2}{1+|t-r|}\sum_{|\beta_2'|=|\beta_2|+1}|Z^{\beta_2'}\p v|.
$$
Because of $|\beta|<|\al|\leq 4$, (2.29), and the fact that
$\left\|(1+|t-r|^{-1}f)(t, \cdot)\right\|_{L^2}\leq\|\p
f(t,\cdot)\|_{L^2}$ for the function $f(t,x)\in C^1(\Bbb R^+\times\Bbb
R^2)$ with $\text{supp}f\subseteq\{r\leq M+t\}$ (this inequality can
be found in [22]), the first term in the right-hand side of (2.35) can
be estimated as
$$
\align
\int_{\Bbb R^2}|Z^{\beta_1}v| \cdot|Z^{\beta_2}\p_i^2 v|\cdot|\p_t
Z^{\al}v|\,dx
& \leq C_b\sum_{|\beta_2'|=|\beta_2|+1}\int_{\Bbb
R^2}\f{1}{1+|t-r|}|Z^{\beta_1}v|\cdot|Z^{\beta'_2}\p v|\cdot|\p_t
Z^{\al}v|\,dx\\
& \leq \f{C_b\ve}{\sqrt{1+t}}\,E(t).\tag2.36
\endalign
$$

Analogously,
$$
\int_{\Bbb
R^2}|Z^{\beta_1}u_a|\cdot|Z^{\beta_2}\p_i^2 v|\cdot|\p_t
Z^{\al}v|\,dx\le\f{C_b\ve}{\sqrt{1+t}}\,E(t).
$$

Therefore, one obtains
$$
\dsize\sum_{|\beta|<|\al|} \int_{\Bbb
R^2}|Z^{\beta}\bigl(c_1u\p_1^2v+c_2u\p_2^2v\bigr)|\cdot|\p_t Z^\alpha
v|\,dx\leq \f{C_b\ve}{\sqrt{1+t}}\,E(t).\tag 2.37
$$

\subhead (B) Estimation of $\ds\int_{\Bbb
R^2}\left|\bigl[Z^{\al},c_1u\p_1^2+c_2u\p_2^2]v\right|\cdot|\p_t
Z^{\al}v|\,dx$ \endsubhead
For $i=1,2$,
$$
\multline \int_{\Bbb
R^2}\left|\bigl[Z^{\al},u\p_i^2]v\right|\cdot|\p_t
Z^{\al}v|\,dx\\
\aligned
& \leq C_b\sum_{\Sb|\al_1|+|\al_2|=|\al|\\|\al_1|\geq
1\endSb}\int_{\Bbb R^2}|Z^{\al_1}u|\cdot|Z^{\al_2}\p_i^2
v|\cdot|\p_t Z^{\al}v|\,dx+\int_{\Bbb R^2}|u|\cdot
\left|[Z^\al,\p_i^2]v\right|\cdot|\p_tZ^{\al}v|\,dx\\
& \leq C_b\left(\sum_{\Sb|\al_1|+|\al_2|=|\al|\\|\al_1|\geq
1\endSb}\int_{\Bbb R^2}|Z^{\al_1}u_a|\cdot|Z^{\al_2}\p_i^2
v|\cdot|\p_t
Z^{\al}v|\,dx+\sum_{|\beta|<|\al|}\int_{\Bbb R^2}|u|\cdot|\p^2Z^\beta v|\cdot
|\p_tZ^{\al}v|\,dx \right.\\
& \qquad \left. +\sum_{\Sb|\al_1|+|\al_2|=|\al|\\|\al_1|\geq 1\endSb}\int_{\Bbb
R^2}|Z^{\al_1}v|\cdot|Z^{\al_2}\p_i^2 v|\cdot|\p_t
Z^{\al}v|\,dx\right).
\endaligned
\endmultline
$$
By the same argument as in (2.37), one then has
$$
\ds\int_{\Bbb
R^2}\left|\bigl[Z^{\al},c_1u\p_1^2+c_2u\p_2^2]v\right|\cdot|\p_t
Z^{\al}v|\,dx\leq \f{C_b\ve}{\sqrt{1+t}}\,E(t).\tag 2.38
$$

\smallskip

Next we treat each of the terms $\ds\int_{\Bbb
R^2}|Z^{\beta}F|\cdot|\p_t Z^{\al}v|\,dx$, $|\beta|\leq |\al|$, which
are included in $\dsize\sum_{|\al|\leq 4}\int_{\Bbb R^2}|G|\cdot|\p_t
Z^\al v|\,dx$.

\subhead (C) Estimation of $\ds\int_{\Bbb
R^2}|Z^{\beta}J_a|\cdot|\p_t Z^{\al}v|\,dx$ \endsubhead
In this case, one has
$$\int_{\Bbb R^2}|Z^{\beta}J_a|\cdot|\p_t Z^{\al}v|\,dx\leq
\|Z^{\beta}J_a\|_{L^{2}}\,\sqrt{E(t)}.\tag 2.39$$

\subhead (D) Estimation of $\ds\int_{\Bbb
R^2}|Z^{\beta}\bigl(c_1v\p_1^2u_a+c_2v\p_2^2u_a\bigr)|\cdot|\p_t Z^{\al}v|\,dx$
\endsubhead
Due to (2.34) a direct computation yields, for $i=1,2$,
$$\align
\int_{\Bbb R^2}|Z^{\beta} \bigl(v\p_i^2u_a\bigr)|\cdot|\p_t Z^{\al}v|\,dx
&\leq C_b\sum_{|\beta_1|+|\beta_2|=|\beta|}\int_{\Bbb
R^2}|Z^{\beta_1}v|\cdot|Z^{\beta_2}\p_i^2 u_a|\cdot|\p_t
Z^{\al}v|\,dx\\
&\leq C_b \sum_{\Sb
|\beta_1|+|\beta_2|=|\beta|\\|\beta_2'|=|\beta_2|+1\endSb}\int_{\Bbb
R^2}\f{1}{1+|t-r|}\,|Z^{\beta_1}v|\cdot|Z^{\beta_2'}\p
u_a|\cdot|\p_t Z^{\al}v|\,dx\tag 2.40 \\
&\leq \f{C_b\ve}{\sqrt{1+t}}\,E(t).
\endalign
$$

\subhead (E) Estimation of $\ds\int_{\Bbb
R^2}|Z^{\beta}\bigl(c_1(\p_1v)^2+c_2(\p_2v)^2\bigr)|\cdot|\p_t
Z^{\al}v|\,dx$ \endsubhead
Similar to (D), one has
$$
\int_{\Bbb R^2}|Z^{\beta}\bigl(c_1(\p_1v)^2+c_2(\p_2v)^2\bigr)|
\cdot|\p_t Z^{\al}v|dx\leq
\f{C_b\ve}{\sqrt{1+t}}\,E(t).\tag 2.41
$$

\subhead (F) Estimation of $\ds\int_{\Bbb
R^2}|Z^{\beta}\bigl(c_1(\p_1 v)(\p_1u_a)
+c_2(\p_2 v)(\p_2u_a)\bigr)|\cdot|\p_t Z^{\al}v|\,dx$ \endsubhead
It follows by direct computation that, for $i=1,2$,
$$\align
\int_{\Bbb R^2}|Z^{\beta}\bigl((\p_i v)(\p_iu_a)\bigr)|\cdot|\p_t Z^{\al}v|\,dx
&\leq C_b\sum_{|\beta_1|+|\beta_2|\leq |\beta|}\int_{\Bbb R^2}|Z^{\beta_1}\p
v|\cdot|Z^{\beta_2}\p u_a|\cdot|\p_t Z^{\al}v|\,dx\\
&\leq \f{C_b\ve}{\sqrt{1+t}}\,E(t).\tag 2.42
\endalign$$

\medskip
Substituting (2.37)-(2.42) into (2.33) yields
$$
E'(t)\leq\f{C_{b}\ve}{\sqrt{1+t}}\,E(t)+\ds\sum_{|\beta|\leq
4}\|Z^\beta J_a(t,\cdot)\|_{L^2}\sqrt{E(t)}.
$$
Thus, by Lemma~2.3 and Gronwall's inequality, one obtains
$$
\|\p Z^\al v(t, \cdot)\|_{L^2}\leq
C_{b}\,\ve^{3/2},\quad|\al|\leq 4,
$$
and further
$$
\|Z^\al\p v(t, \cdot)\|_{L^2}\leq C_{b}\,\ve^{3/2},\quad|\al|\leq
4.\tag2.43
$$

By (2.43) and the Klainerman-Sobolev inequality (see [14,18]), one has
$$
|Z^\kappa\p v|\leq
C_{b}\,\ve^{3/2}(1+t)^{-1/2}(1+|t-r|)^{-1/2},\quad
|\kappa|\leq 2, \enspace t\leq T,\tag 2.44
$$
which means that, for small $\ve$,
$$
|Z^\kappa\p
v|\leq\f{\varepsilon}{2}\,(1+t)^{-1/2}(1+|t-r|)^{-1/2},\quad
|\kappa|\leq 2, \enspace t\leq T.
$$
This completes the proofs of (2.29) and (2.26).
\qed
%%\enddemo

\demo{Proof of the lower bound on $T_{\ve}$}
Lemma 2.4 implies that $\ds\lim_{\overline{\ve\rightarrow
0}}\ve\sqrt{1+ T_\ve}\geq\tau _0$ holds for the lifespan $T_{\ve}$
of smooth solutions to\/ {\rm (1.9)}. Hence,
$$\lim_{\overline{\ve\rightarrow 0}}\ve\sqrt{T_\ve}\geq\tau
_0.\tag2.45$$
which finishes the first part of the proof of Theorem~1.1. \qed
\enddemo

%% ------------------------------------------------------------------------

\head \S3. Proof of Theorem 1.2 \endhead

We will use polar coordinates $(r, \th, t)$ instead of $(x, t)$
to study the problem (1.9) and set
$$
\si=r-t,\quad \tau=\ve\,\sqrt{t}.
$$

Set $u(t,x)=\ds\f{\ve}{\sqrt{r}}\,G(\si,\th,\tau)$ for $r>0$. In this
case, it follows from a direct computation that Eq.~(1.9) takes the
form
$$
\align P(G)\equiv&
-\f{\ve^2}{\sqrt{rt}}\,\p^2_{\si\tau}G-\f{\ve^2}{r}
\left(c_1\cos^2\th+c_2\sin^2\th\right)
G\p_{\si}^2G-\f{\ve^2}{r}\left(c_1\cos^2\th+c_2\sin^2\th\right)
(\p_{\si} G)^2\\
&-\f{\ve^2}{2r^3}\bigl((2c_1-c_2)\cos^2\th+(2c_2-c_1)\sin^2\th\bigr)\,G^2
-\f{\ve^2}{4t^{3/2}r^{1/2}}\,\p_\tau G+\f{\ve^3}{4tr^{1/2}}\,\p_\tau^2G\\
&-\f{\ve^2}{r^3}\left(c_1\sin^2\th+c_2\cos^2\th\right)(\p_\th G)^2
+\f{2\ve^2}{r^2}\sin\th\cos\th(c_1-c_2)(\p_\th G)(\p_\si G)\\
&-\f{\ve}{r^{5/2}}\,\p_\th^2G-\f{\ve}{4r^{5/2}}G-\f{4\ve^2}{r^3}(c_1-c_2)\sin\th
\cos\th G\p_\th G-\f{\ve^2}{r^3}\left(c_1\sin^2\th+c_2\cos^2\th\right)G\p_\th^2G\\
&+\f{\ve^2}{r^2}\left(2c_1+2c_2-3c_1\sin^2\th-3c_2\cos^2\th\right)G\p_\si G
+\f{2\ve^2}{r^2}\sin\th\cos\th(c_1-c_2)G\p_{\si\th}^2G
=0.\tag3.1
\endalign
$$

We introduce an unknown transformation $\Phi$ by
$$
\Phi(s,\th,\tau)
=(\si,\th,\tau),\tag3.2
$$
where $\si=\phi(s,\th,\tau)$, and set
$$
G(\Phi)=v.\tag3.3
$$
Therefore, as $\p_\si G=\p_sv/{\p_s\phi}$, if we can find smooth
functions $\phi$ and $v$ satisfying condition~(H) and (1.15) of
Theorem~1.2, then we will be able to show that the solution $u$ to
(1.9) blows up. Under the transformation (3.2) and (3.3), (3.1) takes
still another form which is explicitly given in the following lemma:

\proclaim{Lemma 3.1} Let $R=1+\ds\f{\ve^2\phi}{\tau^2}$. Then one has
$$
-\f{r}{\ve^2}P(G)\equiv\ds\f{\p_s^2\phi\p_s v}{(\p_s\phi)^3}\,I_0
+\ds\f{1}{(\p_s\phi)^2}\,I_1
+\ds\f{1}{\p_s\phi}\,I_2+I_3=0,\tag3.4
$$
where
$$
\align
I_0&=-(c_1\cos^2\th+c_2\sin^2\th)v-\ds\f{\ve^4}{R^2\tau^4}
(c_1\sin^2\th+c_2\cos^2\th)v(\p_\th\phi)^2
-\f{2\ve^2}{R\tau^2}(c_1-c_2)\sin\th\cos\th v\p_\th\phi\\
&\qquad +R^{1/2}\p_\tau\phi+\ds\f{\ve^2R^{1/2}}{4\tau}(\p_\tau\phi)^2
-\ds\f{\ve^2}{R^{3/2}\tau^3}(\p_\th\phi)^2,\\
I_1&=-\p_s(\p_sv I_0),\\
I_2&=Z_1\p_s v+\ve^2\p_s v N\phi+
\ve^2\p_s vh_1(\ve,\th,\tau,v,\p_{\th}v,\phi,\p_{\th}\phi,\p_{\tau}\phi),\\
I_3&=-\ve^2Nv+\ve^2h_2(\ve,\th,\tau,v,\p_{\th}v,\p_{\tau}v,\phi),\\
\endalign
$$
$h_1,h_2$ are smooth functions the explicit expression of which is
not required, and the first-order differential operator $Z_1$ and the
second-order differential operator $N$, respectively, are of the form
$$
\align
Z_1&=R^{1/2}\bigl(1+\ds\f{\ve^2}{2\tau}\p_\tau\phi\bigr)\p_\tau
-\ds\f{2\ve^2}{R\tau^2}\biggl(\ds\f{\ve^2}{R\tau^2}
(c_1\sin^2\th+c_2\cos^2\th)v\p_\th\phi+\sin\th\cos\th (c_1-c_2)v
+\ds\f{\p_\th\phi}{R^{1/2}\tau}\biggr)\p_\th\\
&\equiv\delta_1\p_\tau+\ve^2\delta_2\p_{\th},\\
N&=\ds\f{R^{1/2}}{4\tau}\p_\tau^2-\ds\f{1}{R^{3/2}\tau^3}\biggl(
1+\ds\f{\ve^2}{R^{1/2}\tau}(c_1\sin^2\th+c_2\cos^2\th)v\biggr)\p_\th^2\\
&\equiv N_1\p_\tau^2+N_2\p_{\th}^2,\\
\endalign
$$
where
$$
\align
\delta_1&=R^{1/2}\bigl(1+\ds\f{\ve^2}{2\tau}\p_\tau\phi\bigr),\quad
\delta_2=-\ds\f{2\ve^2}{R\tau^2}\biggl(\ds\f{\ve^2}{R\tau^2}
(c_1\sin^2\th+c_2\cos^2\th)v\p_\th\phi+\sin\th\cos\th (c_1-c_2)v
+\ds\f{\p_\th\phi}{R^{1/2}\tau}\biggr),\\
N_1&=\ds\f{R^{1/2}}{4\tau},\quad N_2=-\ds\f{1}{R^{3/2}\tau^3}\biggl(
1+\ds\f{\ve^2}{R^{1/2}\tau}(c_1\sin^2\th+c_2\cos^2\th)v\biggr).
\endalign
$$
\endproclaim

It follows from Lemma~3.1 that, in order to solve the nonlinear
equation $P(G)=0$, it suffices to solve the system
$$
\cases
&I_0=0,\\
&I_2+\p_s\phi I_3=0,\\
\endcases\tag3.5
$$
which is also called the blowup system for (1.9) in the terminology of
[4-5] (where nonlinear wave equations such as $\p_t^2v-\Delta_xv
+\ds\sum_{0\le i,j,k\le 2}g_{ij}^k\p_k v\p_{ij}^2v=0$ are dealt with).

\smallskip

The related process is divided into the six parts.

\subhead (A) Local existence of a solution to (3.5) \endsubhead
>From the explicit expression of $I_0$, one has that $\ds\f{\p
I_0}{\p(\p_\tau\phi)}=R^{1/2}+\ds\f{\ve^2R^{1/2}}{2\tau}\p_\tau\phi>0$
for $\ve>0$ small and $\phi$ a smooth function. By the implicit
function theorem, one then obtains from the equation $I_0=0$ that
$$
\p_\tau\phi=E(\ve,\th,\tau,v,\phi,\p_{\th}\phi),\tag3.6
$$
where $E$ is a smooth function of its arguments.

By \S 2, for $C_0>0$ large enough and $\eta>0$ sufficiently small, one
also has that the equation $P(G)=0$ can be solved for $G$ in a strip
$$
D_S=\left\{(\si,\th,\tau)\mid\si\in[-C_0,M],\,
\th\in[\th^0-\dl_0,\th^0+\dl_0],\,\tau\in[\tau_1,\tau_1+\eta]\right\}
$$
with initial data $\ds\f{\sqrt{r}}{\ve}u(t,x)$ given at time
$t=(\tau_1/\ve)^2$ (since (1.9) has a unique smooth solution
there). Here, $\tau_1>0$ is a fixed constant satisfying
$\tau_1<\tau_0$, and $\dl_0>0$ and $0<\eta<\tau_0-\tau_1$ are
sufficiently small.

For $\eta>0$ sufficiently small, Eq.~(3.5) then has a unique solution
$\overline{\phi}$ with initial data $\overline\phi(s,\th,\tau_1)=s$
(note that the smooth solution $u(t,x)$ of (1.9) exists for $t\le
((\tau_1+\eta)/\ve)^2$, as $G(\si,\th,\tau)$ exists for
$\tau\le\tau_1+\eta$).

Setting $\overline{v}=G(\overline{\phi},\th,\tau)$ in the strip $D_S$,
one hence gets a local solution to the blowup system (3.5). Moreover,
from the uniqueness result on the solution $u(t,x)$ to (1.9) for $t\in
\left[0, ((\tau_1+\eta)/\ve)^2\right]$, one has that $\overline{v}$ and
$\overline{\phi}-s$ are smooth and flat on $\{s=M\}$.

\subhead (B) Choice of the domain and the scalar equation for $\phi$
\endsubhead
As in [4-5], in order to obtain a weighted energy estimate on the
linearized system of (3.5) on a suitable domain $D$, we choose a
``nearly horizontal'' surface $\Sigma$ through $\{\tau=\tau_1,s=M\}$
as part of the boundary of $D$, where $\Sigma$ is the characteristic
surface of the operator $Z_1\p_s-\ve^2\p_s{\bar\phi}N$ the
coefficients of which are computed using
$(\overline{v}, \overline{\phi})$. Let $\tau=\psi(s,\th)+\tau_1$ be
the equation of $\Sigma$, where $\psi(M, \th)=0$. Then, in view of
part (A) and for small $\ve>0$, $\nabla^{\al}_{s,\th}\psi=O(\ve^2)$
and $\p_s\psi\leq 0$ holds in $D_S$ for $\al\in {\Bbb N_0^2}$.

We choose a cutoff function $\chi\in C^\infty(\Bbb R)$ with
$\chi(p)=1$ for $p\leq \f{1}{2}$, and $\chi(p)=0$ for $p\geq 1$ and
make the change of variables
$$
X=s,\quad Y=\th,\quad T=\tau-\tau_1-\psi(s,\th)
\chi\left(\f{\tau-\tau_1}{\eta}\right).\tag3.7
$$
The surface $\Sigma$ then becomes $\{T=0\}$. We will work in the
domain $D_1=\left\{(X, Y, T)\mid-C_0\le X\le M,\, \right.$
$\left.\th^0-\th_0\le Y\le\th^0+\th_0,\, 0\le
T\le\tau_\ve-\tau\right\}$. Note that $D_1$ is actually unknown at the
moment, as we do not know the precise value of $\tau_{\ve}$ yet.

Next we derive from (3.5) a scalar equation for $\phi$ in the new
coordinate system (3.7).
Since $\p_v I_0\neq0$ for small $\ve>0$, it follows from $I_0=0$ that
$v$ can be expressed as
$$
v=F(\ve,\th,\tau,\phi,\p_{\th}\phi,\p_\tau\phi),\tag3.8
$$
where $F$ is a smooth function of its arguments.
Substituting (3.8) into the second and third equation of (3.5) and
going through the direct computations yields
$$
L(\phi)\equiv\widetilde{Z_1}S F-\ve^2(S\phi)
\widetilde{N}F+\ve^2(SF)\widetilde{N}\phi+
\ve^2(SF)\widetilde{h_1}+\ve^2(S\phi)\widetilde{h_2}=0,
\tag3.9
$$
where
$$
\align
\widetilde{Z_1}&=\widetilde{\delta}_1\p_T+\ve^2\widetilde{\delta}_2\p_{Y},
\quad
\widetilde{N}=\widetilde{N}_1\p_T^2+2\ve^2\widetilde{N}_2\p_{YT}^2
+\widetilde{N}_3\p_{Y}^2,\\
\widetilde{h}_1&=h_1-N_1\psi\chi''\p_T\phi/{\eta^2}-N_2\p_{\th}^2
\psi\chi\p_T\phi,\\
\widetilde{h}_2&=h_2+N_1\psi\chi''\p_TF/{\eta^2}+N_2\p_{\th_1}^2
\psi\chi\p_TF,\\
S&=\p_s=\p_X-\p_s\psi\chi\p_T, \\
\intertext{where}
\widetilde{\delta}_1&=\delta_1\p_\tau T+\ve^2\delta_2\p_{\th} T,
\quad \widetilde{\delta}_2=\delta_2,\\
\widetilde{N}_1&=N_1(\p_\tau T)^2+N_2(\p_{\th}T)^2,\quad
\widetilde{N}_2=\ve^{-2}N_2\p_{\th}T,\quad
\widetilde{N}_3=N_2.
\endalign
$$

In order to solve the blowup system (3.5), one hence only needs to
solve (3.9) because of (3.8). As in [4-5], we will use the
Nash-Moser-H\"ormander iteration method to solve Eq.~(3.9) under the
restriction (H) of Theorem 1.2.

\subhead (C) The construction of an approximate solution to (3.9)
and the condition (H) \endsubhead
As a first step to use the Nash-Moser-H\"ormander iteration method,
one needs to construct an approximate solution $\phi_a$ to (3.9) such
that $\phi_a$ satisfies (H) of Theorem 1.2 near some point
$m_\epsilon$.

For $\ve=0$, the blowup system (3.5) becomes
$$
(c_1\cos^2Y+c_2\sin^2Y)v=\p_T\phi,\quad \p_Tv=0\tag3.10
$$
with the initial value conditions
$$
\phi(X,Y,0)=X,\quad v(X,Y,0)=F_0(\si(X,Y,\tau_1),Y)\tag3.11
$$
and the boundary condition
$$v|_{X=M}=0,\tag3.12$$
where the function $\si(X,Y,\tau_1)$ in (3.11) is determined by
$X=\si+F_0(\si,Y)\tau_1(c_1\cos^2Y+c_2\sin^2Y)$.

>From (3.10)-(3.12), one finds a solution to (3.9) for $\ve=0$, namely
$$
\overline{\phi}_0(X,Y,T)=X+T(c_1\cos^2Y+c_2\sin^2Y)F_0(\si(X,Y,\tau_1),Y).
\tag3.13
$$

Note that (3.9) admits a local solution $\overline{\phi}$ for $0\le
T\le \eta$ the existence of which has been proven in part~(A). Upon
glueing $\overline{\phi}$ and $\overline{\phi}_0$ one obtains an
approximate solution to (3.9), namely
$$
\phi_a(X,Y,T)=\chi\left(\f{T}{\eta}\right)\overline{\phi}(X,Y,T)
+\biggl(1-\chi\left(\f{T}{\eta}\right)\biggr)
\overline{\phi}_0(X,Y,T).\tag3.14
$$

By a direct verification, one has $L(\phi_a)=f_a$, where $f_a$ is smooth,
flat on $\{X=M\}$, and zero near $\{T=0\}$.

In addition, under the assumption (1.11), one can
show that $\phi_a$ satisfies (H) at the point $(\bar\si^0, \th^0,
\tau_0-\tau_1)$ with $\bar\si^0=\si^0
+(c_1\sin^2\th^0+c_2\cos^2\th^0)F_0(\si^0,
\th^0)\tau_1$:

\proclaim{Lemma 3.2} The approximate solution $\phi_a$ constructed in
\rom{(3.14)} satisfies \rom{(H)} near the point $(\bar\si^0, \th^0,
\tau_0-\tau_1)$.
\endproclaim

\demo{Proof} Note that
$$
\overline{\phi}_0(X,Y,T)=X+T(c_1\cos^2Y+c_2\sin^2Y)F_0(\si(X,Y,\tau_1),Y),
$$
where $\si(X,Y,\tau_1)$ is determined from the expression
$X=\si+F_0(\si,Y)\tau_1(c_1\cos^2Y+c_2\sin^2Y)$ (this follows as in the
proof of Lemma 2.1).

Set $W(X,Y)=F_0(\si(X,Y,\tau_1), Y)(c_1\cos^2\th+c_2\sin^2\th)$.
First, we assert that
$$
\p_XW(\bar\si^0,\th^0)=\min \p_XW(X,Y).\tag3.15
$$
Indeed, it follows from (1.11) and a direct computation that
$$
\align
\nabla_{X,Y}\p_XW(\bar\si^0, \th^0)&=0, \\
\nabla_X^2\p_XW(\bar\si^0, \th^0)&=
\ds\f{(c_1\cos^2\th^0+c_2\sin^2\th^0)\p_{\si}^3F_0(\si^0,\th^0)}
{(1+\tau_1(c_1\cos^2\th^0+c_2\sin^2\th^0)\p_{\si}F_0(\si^0,\th^0))^4}, \\
\nabla_{XY}^2\p_XW(\bar\si^0, \th^0) &=
\ds\f{-\tau_1\p_{\th}((c_1\cos^2\th+c_2\sin^2\th)F_0)(\si^0,\th^0)}
{(1+\tau_1(c_1\cos^2\th^0+c_2\sin^2\th^0)\p_{\si}F_0(\si^0,\th^0))^4}
\p_{\si}^3F_0(\si^0,\th^0)(c_1\cos^2\th^0+c_2\sin^2\th^0) \\
& \qquad
+\ds\f{\p_{\th}\big(\p_{\si}^2F_0(c_1\cos^2\th+c_2\sin^2\th)\big)(\si^0,\th^0)}
{(1+\tau_1(c_1\cos^2\th^0+c_2\sin^2\th^0)\p_{\si}F_0(\si^0,\th^0))^3},
\endalign
$$
and
$$
\align
\nabla_{Y}^2\p_XW(\bar\si^0, \th^0)&=
\ds\f{\tau_1^2\big(\p_{\th}((c_1\cos^2\th+c_2\sin^2\th)F_0)(\si^0,\th^0)\big)^2}
{(1+\tau_1(c_1\cos^2\th^0+c_2\sin^2\th^0)\p_{\si}F_0(\si^0,\th^0))^4}
(c_1\cos^2\th^0+c_2\sin^2\th^0)\p_{\si}^3F_0(\si^0,\th^0) \\
& \qquad -\ds\f{2\tau_1\p_{\th}((c_1\cos^2\th+c_2\sin^2\th)F_0)(\si^0,\th^0)}
{(1+\tau_1(c_1\cos^2\th^0+c_2\sin^2\th^0)\p_{\si}F_0(\si^0,\th^0))^3}
\p_{\th}((c_1\cos^2\th+c_2\sin^2\th)\p_\si^2F_0)(\si^0,\th^0) \\
& \qquad +\ds\f{\p_{\th}^2((c_1\cos^2\th+c_2\sin^2\th)\p_{\si}F_0)(\si^0,\th^0)}
{(1+\tau_1(c_1\cos^2\th^0+c_2\sin^2\th^0)\p_\si F_0(\si^0,\th^0))^{2}}.
\endalign
$$
This, together with $\na_{\si, \th}^2[\p_{\si}F_0(\si,
\th)(c_1\cos^2\th+c_2\sin^2\th)]|_{(\si, \th)=(\si^0, \th^0)}>0$, yields by a direct,
but tedious computation
$$
\nabla_{X,Y}^2\p_XW(\bar\si^0,\th^0)>0.\tag 3.16
$$
Thus, the assertion (3.15) has been shown. Moreover, by the uniqueness
of the minimum point of the function
$\p_{\si}F_0(\si,\th)(c_1\cos^2\th+c_2\sin^2\th)$, one has that
$(\bar\si^0,\th^0,\tau_1)$ is also the unique minimum point of
$\p_XW(X,Y)$.

\smallskip

We now establish that $\phi_a$ satisfies (H) near the point
$(\bar\si^0,\th^0,\tau_0-\tau_1)$.

\roster
\item"(i)"  By $\p_X\overline{\phi}(X,Y,0)=1$ and the smallness of $\eta>0$,
one can assume that, for $T\leq\eta$,
$$
\p_X\overline{\phi}(X,Y,T)>0.
$$
In addition,
$$
\p_X\overline{\phi}_0(X,Y,T)=1+T\p_XW(X,Y).\tag 3.17
$$
If $\p_XW(X,Y)\geq 0$ at some point $(X,Y,\tau_1)$, then
$\p_X\overline{\phi}_0(X,Y,T)\geq 1$.
If $\p_XW(X,Y)<0$ at some point $(X,Y,\tau_1)$, then due to the fact
$\ds T\leq T_0= -\min \frac1{\p_X \big(F_0(\si(X,Y,\tau_1)\big)
\left(c_1\cos^2Y+c_2\sin^2Y\right)}=\tau_0-\tau_1$ one has from (3.17)
that
$$
\p_X\overline{\phi}_0(X,Y,T)\geq 1+\p_XW(X,Y)T_0\geq
0.
$$
Consequently, $\p_X\phi_a(X,Y,T)\geq 0$ holds.

On the other hand, $\p_X\phi_a(X,Y,T)=0$ holds if and only if
$T\geq\eta$ and $\p_X\overline{\phi}_0(X,Y,T)$ $=0$ which gives
$$
\p_X\phi_a(X,Y,T)=0\enspace \Longleftrightarrow \enspace (X,Y,T)
=(\bar\si^0,\th^0,\tau_0-\tau_1).
$$

\item"(ii)" It follows from the expression for $\phi_a$ and the smallness of
$\eta>0$ that in the neighborhood of $(\bar\si^0,\th^0,\tau_0-\tau_1)$
$$
\phi_a(X,Y,T)=X+TW(X,Y)\tag 3.18
$$
which gives $\p_{XT}^2\phi_a(\bar\si^0,\th^0,\tau_0-\tau_1)<0$.
In addition, in view of $\nabla_{X,Y}\p_XW(X,Y)(\bar\si^0, \th^0)=0$,
(3.16), and (3.18), one readily obtains
$$
\nabla_{X,Y}\p_X\phi_a(\bar\si^0,\th^0,\tau_0-\tau_1)=0,
\qquad\nabla^2_{X,Y}\p_X\phi_a(\bar\si^0,\th^0,\tau_0-\tau_1)>0.
$$
\endroster

Collecting all the assertions above concludes the proof of
Lemma~3.2. \qed
\enddemo

\subhead (D) Goursat problem for the nonlinear equation (3.9)  on a
fixed domain \endsubhead
In order to adjust the height of the domain $D_1$ as in [4] we perform
a change of variables depending on a parameter $\la$ close to zero,
$$
X=x,\quad Y=y,\quad
T=T(\rho,\la)=(\tau_0-\tau_1)(\rho+\la\rho(1-\chi_1(\rho))),\tag3.19
$$
where $\chi_1$ is $1$ near $0$ and $0$ near $1$.
>From now on we will be working on a fixed subdomain of $D_1$,
$$
D_2=\left\{(x,y,\rho)\mid -C_0\le x\le M,\,
\th^0-\dl_0\le y\le\th^0+\dl_0,\, 0\le\rho\le 1\right\}
$$
and write Eq.~(3.9) as
$$
L(\la,\phi)=0.\tag3.20
$$

For $\la=\la_0=0$, the approximate solution to (3.20) is
$$
\phi_0(x,y,\rho)=\phi_a(x,y,T(\rho,0))=\phi_a(x,y,(\tau_0-\tau_1)\rho),
$$
where $L(\la_0,\phi_0)=f_{0}(x,y,\rho)=
f_{a}(x,y,(\tau_0-\tau_1)\rho)$. Moreover, $\phi_0$ satisfies (H) in
$D_2$ at some point $(x_0,y_0,1)$ by part~(C).

On the characteristic surfaces $\{x=M\}$ and $\{\rho=0\}$ of
Eq.~(3.20), we impose the natural boundary conditions
$$
\text{$\phi$ is flat on $\{x=M\}$ and $\phi-\phi_0$ is flat on
$\{\rho=0\}$, respectively.} \tag3.21
$$

\subhead (E) Linearizing (3.20) under the condition (H) \endsubhead
In order to solve (3.20) together with (3.21) in the domain $D_2$
under the condition (H), we are required to linearize (3.20) suitably.

Denote the linearized operator of $L$ by
$$
L'(\la,\phi)(\dot{\la},\dot{\phi})=\p_\la
L(\la,\phi)\dot{\la}+\p_\phi
L(\la,\phi)\dot{\phi}=\dot{f}.\tag3.22
$$
In addition, if $L(\la,\phi)=f$, then taking the derivative with
respect to the variable $\la$ yields
$$
\p_\la L(\la,\phi)+\p_\phi L(\la,\phi)\biggl(\p_{\rho}\phi\f{\p_\la
T}{\p_{\rho} T}\biggr)=\p_{\rho}f\f{\p_\la T}{\p_{\rho}
T}.\tag3.23
$$
Therefore, if one wants to solve $L(\la, \phi)=f$ for a small
right-hand side $f$, then it follows from the standard
Nash-Moser-H\"ormander iteration method that we are only required to
solve the linearized equation $L'(\la,\phi)(\dot{\la},\dot{\phi})=\dot
f$ and provide the needed tame estimate (see [3]). From (3.22)-(3.23),
one has $L'(\la,\phi)(\dot{\la},\dot{\phi})=\p_\phi
L(\la,\phi)\biggl(\dot\phi-\dot\la
\p_{\rho}\phi\ds\f{\p_\la T}{\p_{\rho} T}\biggr)+\dot\la
\p_{\rho}f\ds\f{\p_{\la}T}{\p_{\rho} T}$. Setting
$\dot\Phi=\dot\phi-\dot\la \p_{\rho}\phi\ds\f{\p_\la T}{\p_{\rho}
T}$, it suffices to solve the equation
$$
\cases
&\p_\phi L(\la,\phi)\dot\Phi=\dot f,\\
&\text{$\dot\Phi$ is flat on both $\{x=M\}$ and $\{\rho=0\}$}
\endcases\tag 3.24
$$
for a right-hand side $\dot f$ which is also flat on both $\{x=M\}$
and $\{\rho=0\}$, since the second-order error term (here $\dot\la
\p_{\rho}f\ds\f{\p_{\la}T}{\p_{\rho} T}$) does not play an
essential role in the Nash-Moser-H\"ormander iteration (see [3]).

It follows from a direct, but tedious computation concerning $\p_\phi
L_i(\la,\phi)$ that from (3.24) one obtains
$$
\cases
&\overline{P}\dot\Phi\equiv ZSZ\dot\Phi-\ve^2(S\phi)QZ\dot\Phi
+\ve^2l(\dot\Phi)=\dot {f},\\
&\text{$\dot\Phi$ is flat on both $\{x=M\}$ and $\{\rho=0\}$}
\endcases\tag3.25
$$
as the linearized problem of (3.20), where
$$
Z=\p_\rho+\ve^2z_0\p_y, \quad S=\p_x+\ve^2s_0\p_\rho, \quad
Q=Q_1Z^2+2\ve^2Q_2Z\p_{y}+Q_3\p_{y}^2,
$$
here $z_0, s_0$, and $Q_i$ are smooth. More specifically,
$$
\align
z_0&=z_0(x,\rho,\la,\phi,\p_{y}\phi,\p_\rho\phi),\quad
s_0=s_0(x,y,\rho,\la),\\
Q_1&=\ds\f{1}{4(\tau_1+T)\p_\rho T}+O(\ve^2),\quad
Q_3=-\ds\f{\p_\rho T}{(\tau_1+T)^3}+O(\ve^2),
\endalign
$$
and $l$ is a second-order operator which is a linear combination of
$id,S,Z,\p_{y},SZ,Z^2,Z\p_{y},\p_{y}^2$ and whose coefficients depend
on the derivatives of $\phi$ up to third order.

\subhead (F) The tame estimate and
solvability of (3.24) \endsubhead
Comparing the operator $\bar P$ in (3.25) with the operator
$\p_{\phi}{\Cal L}(\la, \phi)$ in Proposition IV.1 of [4], one sees
that $\bar P$ is just of the form of $\p_{\phi}{\Cal L}(\la, \phi)$
with $B\equiv 0$ and $b_0\equiv 0$. By carefully checking the proofs
of Proposition~IV.2.2, Proposition~IV.3.1, and Proposition~IV.4 of
[4], one then has under the condition (H) on the function $\phi$ near
some point $(\bar x_0,\bar y, 1)$:

\proclaim{Lemma 3.3} There exists a subdomain $D_0$ of $D_2$ which is a
domain of influence domain for the first-order differential operator
$\t Z_1$ in \rom{(3.9)}, that contains the point $(\bar x_0, \bar y,
1)$, and that is bounded by the planes $\{x=-C_0\}$, $\{x=M\}$,
$\{\rho=0\}$, $\{\rho=1\}$ with the followoing property\rom: If
$|\phi-\phi_0|_{C^7(D_3)}\leq\ve_0$ with $\ve_0$ a small positive
constant and if $f\in C^\infty(D_3)$ is flat on both $\{x=M\}$ and
$\{\rho=0\}$, then \rom{(3.25)} has a unique smooth solution in
$D_3$. Moreover, one has the energy estimate \rom(i.e., tame
estimate\rom)
$$
|\dot\Phi|_s\leq C_s\left(|\dot f|_{s+n_0}
+|\dot f|_{n_0}(1+|\phi|_{s+n_0})\right),\tag3.26
$$
for any $s\in\Bbb N$, where $|\cdot|_s=\|\cdot\|_{H^s(D_3)}$ and
$n_0\in\Bbb N$ is some fixed integer.
\endproclaim

Based on Lemma~3.3 and the standard Nash-Moser-H\"ormander
iteration method (see [3, 4-5]), and using
$\p_Xv(\bar\si^0, \th^0)\neq0$ in (3.11) and the Sobolev imbedding
theorem, we have now completed the proof of Theorem~1.2.

%% ------------------------------------------------------------------------

\head \S4. Proof of Theorem 1.1 \endhead

Using Theorem 1.2, we now conclude the proof of Theorem 1.1.

Recall that so far we have obtained the $C^3$ solution $\phi$ to
Eq.~(3.20), (3.21) in the domain $D_0$. By (3.8), we immediately
obtain $v$ in $D_0$. Indeed, we have solved the modified blowup system
(3.5) in $D_0$. Therefore, the solution to (1.9) is obtained in the
domain $\Phi(D_0)$ in the coordinate system $(s,\th,\tau)$, hence
$\|u\|_{C(\Phi(D_0))}\leq C\ve^2$. Now we go back to the original
coordinate system $(r, \th, t)$ so that the conclusions of Theorem~1.1
can be obtained.

For $(s,\th,\tau)\in D_0$ close to the point $m_\ve$ given in
Theorem~1.2, by Taylor's formula and condition (H) in Theorem~1.2,
there exists a point $(\bar s, \bar\th, \bar\tau) =(\bar\la
s+(1-\bar\la)s_{\ve}, \bar\la \th+(1-\bar\la)\th_{\ve},
\bar\la \tau+(1-\bar\la)\tau_{\ve})$ with $0<\bar\la<1$
such that
$$
\multline
\p_s\phi(s,\th,\tau)=\p_{s\tau}^2\phi(m_{\ve})(\tau-\tau_\ve)
\\
+\ds\f12 (s-s_\ve, \th-\th_\ve, \tau-\tau_\ve)\na^2_{s,\th,\tau}\p_s
\phi(\bar s, \bar\th, \bar\tau)(s-s_\ve, \th-\th_\ve, \tau-\tau_\ve)^T.
\qquad \text{$ $}
\endmultline
\tag4.1
$$
In addition, we may assume $-2c_0\leq\p_{s\tau}^2\phi\leq-c_0$ in
$D_0$ since $\p_{s\tau}^2\phi(m_{\ve})<0$ and $\phi\in C^3(D_0)$, here
$c_0>0$ is a constant. Together with
$\na_{s,\th}^2\p_s\phi(m_{\ve})>0$, for $(s,\th,\tau)\in D_0$, this
yields
$$
\p_s\phi(s,\th,\tau)\geq
c_0(\tau_\ve-\tau)=c_0\ve\,\ds\f{T_\ve-t}{\sqrt{T_\ve}
+\sqrt{t}}\geq\f{c_0\ve}{4}\cdot\f{T_\ve-t}{\sqrt{t}}.\tag4.2
$$

Furthermore, if $|(s-s_\ve,\th-\th_\ve)|<\tau_{\ve}-\tau$, then
$$
|\p_s\phi(s,\th,\tau)|\le 3c_0(\tau_{\ve}-\tau)
\leq c_0\ve\,\f{T_\ve-t}{\sqrt t}.\tag4.3
$$
>From the expression for $u=\ds\f{\ve}{\sqrt r}\,v$, one has
$$
\align
&\p_1u=-\ds\f{\ve\cos\th}{2r^{3/2}}\,v-\f{\ve}{\sqrt r}\left(\p_\th v
-\f{\p_sv}{\p_s\phi}\,\p_\th\phi\right)\f{\sin\th}{r},\\
&\p_2u=-\ds\f{\ve\sin\th}{2r^{3/2}}\,v+\f{\ve}{\sqrt r}\left(\p_\th v
-\f{\p_sv}{\p_s\phi}\,\p_\th\phi\right)\f{\cos\th}{r},\\
&\p_tu=\f{\ve^2}{2\sqrt{rt}}\left(\p_\tau v-\f{\p_sv}{\p_s\phi}\,
\p_\tau\phi\right).
\endalign
$$
Substituting (4.2)-(4.3) into the these formulas yields
$$
\f{1}{C\left(T_\ve-t\right)}\leq\|\partial_t u\|_{L^\infty(\Phi(D_0)}
\enspace\text{and}\enspace
\|\nabla_{t,x} u(t,\cdot)\|_{L^\infty(\Phi(D_0))}
\leq\f{C}{T_\ve-t}.\tag4.4
$$

Owing to assumption (1.11), outside $\Phi(D_0)$ and for $t\le
T_{\ve}$, the smooth solution of (2.1) does not blow up in $(\{t\le
T_{\ve}\}\times\Bbb R^3)\setminus \Phi(D_0)$. Therefore, similar to
the proof of Lemma 2.4, one obtains in $(\{t\le T_{\ve}\}\times\Bbb
R^3)\setminus \Phi(D_0)$ that
$$
\text{$|\p u|\leq C\ve\,(1+t)^{-1/2}$ and $|u|\leq C\ve$.}
$$

Finally, by Theorem 1.2 and the related Nash-Moser-H\"ormander
iteration process, one concludes that $\ds\lim_{\ve\to
0}\tau_{\ve}=\tau_0$ for the solution $u(t,x)$ when the variables
$(r-t, \th, \ve\sqrt t)$ lie in $\Phi(D_0)$. This implies that the
lifespan $T_{\ve}$ satisfies
$$
\ds\overline{\lim_{\ve\to 0}}\ve \sqrt T_{\ve}\le\tau_0.\tag4.5
$$

Together with (2.45), this yields
$$\ds\lim_{\ve\to 0}\ve \sqrt T_{\ve}=\tau_0,$$
which completes the proof of Theorem 1.1. %%\qed

%% ------------------------------------------------------------------------

%% ------------------------------------------------------------------------


\Refs \refstyle{C}

\ref\key 1 \by R.K.~Agarwal, D.W.~Halt \paper A modified CUSP scheme
in wave/particle split form for unstructured grid Euler flows \inbook
Frontiers of Computational Fluid Dynamics 1994 \eds D.A.~Caughey,
M.M.~Hafez \yr 1995
\endref

\ref\key 2 \by G.~Ali, J.K.~Hunter \paper Diffractive nonlinear geometrical
optics for variational wave equations and the Einstein equations \jour
Comm. Pure Appl. Math. \vol 60 \issue 10 \pages 1522--1557 \yr 2007
\endref

\ref\key 3 \by S.~Alinhac, P.~G\'erard \book Pseudo-differential operators
and the Nash-Moser theorem \transl Translated from the 1991 French
original by S.~Wilson. \bookinfo Grad. Stud. Math. \vol 82 \publ
Amer. Math. Soc. \publaddr Providence, RI \yr 2007
\endref

\ref\key 4 \by S.~Alinhac \paper Blow up of small data solutions for a
class of quasilinear wave equations in two space dimensions \jour
Ann. of Math. (2) \vol 149 \issue 1 \pages 97--127 \yr 1999
\endref

\ref\key 5 \by S.~Alinhac \paper Blow up of small data solutions for a
class of quasilinear wave equations in two space dimensions. II \jour
Acta Math. \vol 182 \issue 1 \pages 1--23 \yr 1999
\endref

\ref\key 6 \by S.~Alinhac \paper An example of blowup at infinity for
quasilinear wave equations \jour Ast\'erisque \vol 284 \pages
1--91 \yr 2003
\endref

\ref\key 7 \by A.~Bressan, Zheng Yuxi \paper Conservative solutions to
a nonlinear variational wave equation \jour Comm. Math. Phys. \vol 266
\issue 2 \pages 471--497 \yr 2006
\endref

\ref\key 8 \by D.~Christodoulou \paper Global solutions of nonlinear
hyperbolic equations for small initial data \jour Comm. Pure Appl.
Math. \vol 39 \issue 2 \pages 267--282 \yr 1986
\endref

\ref\key 9 \by Ding Bingbing, Yin Huicheng \paper On the blowup of
classical solutions to the 3-D pressure-gradient systems \jour
J. Differential Equations \vol 252 \pages 3608--3629 \yr 2012
\endref

\ref\key 10 \by R.T.~Glassey, J.K.~Hunter, Zheng Yuxi \paper
Singularities of a variational wave equation \jour J. Differential
Equations \vol 129 \pages 49--78 \yr 1996
\endref

\ref\key 11 \by H.~Holden, K.H.~Karlsen, N.H.~Risebro \paper A convergent
finite-difference method for a nonlinear variational wave
equation \jour IMA J. Numer. Anal. \vol 29 \issue 3 \pages
539--572 \yr 2009
\endref

\ref\key 12 \by H.~Holden, X.~Raynaud \paper Global semigroup of
conservative solutions of the nonlinear variational wave
equation \jour Arch. Ration. Mech. Anal. \vol 201 \issue 3 \pages
871--964 \yr 2011
\endref

\ref\key 13 \by L.~H\"{o}rmander \paper The lifespan of classical
solutions of nonlinear hyperbolic equations \moreref Mittag-Leffler
report no.~5 \yr 1985
\endref

\ref\key 14 \by L.~H\"{o}rmander \book Lectures on nonlinear hyperbolic
equations \bookinfo Math. Appl. \vol 26 \publ Springer \publaddr
Berlin \yr 1997
\endref

\ref\key 15 \by J.K.~Hunter, R.A.~Saxton\paper Dynamics of director
fields \jour SIAM J. Appl. Math. \vol 51 \pages 1498--1521 \yr 1991
\endref

\ref\key 16 \by J.K.~Hunter, Zheng Yuxi \paper On a nonlinear
hyperbolic variational equation. I. Global existence of weak
solutions \jour Arch. Rational Mech. Anal. \vol 129 \issue 4 \pages
305--353 \yr 1995
\endref

\ref\key 17 \by F.~John \paper Blow-up of radial solutions of $u_{tt}
= c^2(u_t)\Delta u$ in three space dimensions \jour Mat. Apl.
Comput. \vol 4 \issue 1 \pages 3--18 \yr 1985
\endref

\ref\key 18\by S.~Klainerman \paper Remarks on the global Sobolev
inequalities in the Minkowski space $\Bbb R^{n+1}$\jour Comm. Pure
Appl. Math. \vol 40 \pages 111--117 \yr 1987
\endref

\ref\key 19 \by Lei Zhen, Zheng Yuxi \paper A complete global
solution to the pressure gradient equation \jour J. Differential
Equations \vol 236 \pages 280--292 \yr 2007
\endref

\ref\key 20\by Li Fengbai, Wei Xiao \paper Interaction of four
rarefaction waves in the bi-symmetric class of the pressure-gradient
system \jour J. Differential Equations \vol 252 \pages 3920--3952 \yr
2012
\endref

\ref\key 21 \by Li Jun, I.~Witt, Yin Huicheng \paper On the blowup
and lifespan of smooth solutions to a class of 2-D nonlinear wave
equations with small symmetric initial data \finalinfo Preprint \yr
2011
\endref

\ref\key 22 \by H.~Lindblad \paper On the lifespan of solutions of
nonlinear wave equations with small initial data \jour Comm. Pure
Appl.  Math. \vol 43 \issue 4 \pages 445--472 \yr 1990
\endref

\ref\key 23 \by H.~Lindblad \paper Global solutions of quasilinear
wave equations \jour Amer. J. Math. \vol 130 \issue 1 \pages
115--157 \yr 2008
\endref

\ref\key 24 \by R.A.~Saxton \paper Dynamic instability of the liquid
crystal director \inbook Current progress in hyperbolic systems:
Riemann problems and computations (Brunswick, ME, 1988) \pages
325--330 \bookinfo Contemp. Math. \vol 100 \publ Amer. Math. Soc.
\publaddr Providence, RI \yr 1989\endref

\ref\key 25 \by  Song Kyungwoo, Zheng Yuxi \paper Semi-hyperbolic
patches of solutions of the pressure gradient system \jour Discrete
Contin. Dyn. Syst. \vol 24 \issue 4 \pages 1365--1380 \yr 2009
\endref

\ref\key 26 \by Yin Huicheng \paper Formation and construction of a
shock wave for 3-D compressible Euler equations with the spherical
initial data \jour Nagoya Math. J. \vol 175 \pages 125--164 \yr 2004
\endref

\ref\key 27 \by Zhang Ping, Zheng Yuxi \paper Weak solutions to a
nonlinear variational wave equation \jour
Arch. Ration. Mech. Anal. \vol 166 \pages 303--319 \yr 2003
\endref

\ref\key 28 \by Zhang Ping, Zheng Yuxi \paper Conservative solutions
to a system of variational wave equations of nematic liquid
crystals \jour Arch. Ration. Mech. Anal. \vol 195 \issue 3 \pages
701--727 \yr 2010
\endref

\ref\key 29 \by Zheng Yuxi \book Systems of conservation laws:
Two-dimensional Riemann problems \bookinfo Progr. Nonlinear
Differential Equations Appl. \vol 38 \publ Birkh\"auser \publaddr
Boston \yr 2001
\endref

\ref\key 30 \by Zheng Yuxi \paper Two-dimensional regular shock
reflection for the pressure gradient system of conservation laws \jour
Acta Math. Appl. Sin \vol 22 \issue 2 \pages 177--210 \yr 2006
\endref

\ref\key 31 \by Zheng Yuxi, R.~Zachary \paper The pressure gradient
system \jour Methods Appl. Anal. \vol 17 \issue 3 \pages 263--278 \yr
2010
\endref

\endRefs
\enddocument